\documentclass[10pt]{article}

\usepackage{babel}
\usepackage[T1]{fontenc}
\usepackage{epsfig,amsmath,amssymb,amsthm,amsfonts}
\usepackage{lineno}
\usepackage[latin1]{inputenc}
                       
\newtheorem{theorem}{Theorem}
\newtheorem{corollary}{Corollary}

\newtheorem{example}{\bf Example}

\begin{document}

\newcommand{\ncm}{\newcommand}
\ncm{\bfm}[1]{\mbox{\boldmath $#1$}}
\ncm{\sbfm}[1]{\mbox{\scriptsize\boldmath $#1$}}
\ncm{\scr}[1]{\mbox{\scriptsize #1}}
\ncm{\scrmath}[1]{\mbox{\scriptsize $#1$}}
\ncm{\bfmscr}[1]{\mbox{\scriptsize{\boldmath $#1$}}}

\ncm{\R}{{\mathbb{R}}}
\ncm{\Z}{{\mathbb{Z}}}
\ncm{\T}{{\mathbb{T}}}
\ncm{\Smath}{{\mathbb{S}}}
\ncm{\N}{{\mathbb{N}}}
\ncm{\C}{{\mathbb{C}}}
\ncm{\A}{{\mathbb{A}}}
\ncm{\amath}{\bfm{a}}
\ncm{\V}{{\mathbb{V}}}
\ncm{\Hap}{{\mathbb{H}}}
\ncm{\MMD}{{\mathbb{MD}}}

\ncm{\cA}{{\cal A}}
\ncm{\cB}{{\cal B}}
\ncm{\cC}{{\cal C}}
\ncm{\calF}{{\cal F}}
\ncm{\cD}{{\cal D}}
\ncm{\cG}{{\cal G}}
\ncm{\cL}{{\cal L}}
\ncm{\cN}{{\cal N}}
\ncm{\cI}{{\cal I}}
\ncm{\cJ}{{\cal J}}
\ncm{\cH}{{\cal H}}
\ncm{\cV}{{\cal V}}
\ncm{\cW}{{\cal W}}
\ncm{\cT}{{\cal T}}
\ncm{\cX}{{\cal X}}
\ncm{\cQ}{{\cal Q}}
\ncm{\cR}{{\cal R}}
\ncm{\cS}{{\cal S}}
\ncm{\cM}{{\cal M}}
\ncm{\cU}{{\cal U}}
\ncm{\cP}{{\cal P}}
\ncm{\cZ}{{\cal Z}}
\ncm{\cO}{{\cal O}}
\ncm{\cPzer}{{\cal P}_0}
\ncm{\cPone}{{\cal P}_1}
\ncm{\cPk}{{\cP_{\mbox{\scr{known}}}}}
\ncm{\cF}{{\cal F}}
\ncm{\cE}{{\cal E}}
\ncm{\cMD}{{\cal MD}}
\ncm{\tcV}{\tilde{\cal V}}
\ncm{\cCobs}{{\cal C}_{\scr{obs}}}

\ncm{\Om}{\Omega}
\ncm{\om}{\omega}
\ncm{\va}{\varepsilon}
\ncm{\vam}{\varepsilon_{\scr{max}}}
\ncm{\de}{\delta}
\ncm{\De}{\Delta}
\ncm{\ga}{\gamma}
\ncm{\Ga}{\Gamma}
\ncm{\la}{\lambda}
\ncm{\ka}{\kappa}
\ncm{\si}{\sigma}
\ncm{\Si}{\Sigma}
\ncm{\La}{\Lambda}
\ncm{\eps}{\epsilon}

\ncm{\bY}{\bfm{Y}}
\ncm{\bA}{\bfm{A}}
\ncm{\bB}{\bfm{B}}
\ncm{\bC}{\bfm{C}}
\ncm{\bD}{\bfm{D}}
\ncm{\bF}{\bfm{F}}
\ncm{\bI}{\bfm{I}}
\ncm{\bZ}{\bfm{Z}}
\ncm{\bG}{\bfm{G}}
\ncm{\bH}{\bfm{H}}
\ncm{\bL}{\bfm{L}}
\ncm{\bP}{\bfm{P}}
\ncm{\bQ}{\bfm{Q}}
\ncm{\bS}{\bfm{S}}
\ncm{\bT}{\bfm{T}}
\ncm{\bU}{\bfm{U}}
\ncm{\bM}{\bfm{M}}
\ncm{\bN}{\bfm{N}}
\ncm{\bW}{\bfm{W}}
\ncm{\bX}{\bfm{X}}
\ncm{\bu}{\bfm{u}}
\ncm{\bv}{\bfm{v}}
\ncm{\bw}{\bfm{w}}
\ncm{\bwpr}{\bfm{w}^\prime}
\ncm{\bhp}{\bfm{h}^\prime}
\ncm{\bc}{\bfm{c}}
\ncm{\bd}{\bfm{d}}
\ncm{\bh}{\bfm{h}}
\ncm{\bn}{\bfm{n}}
\ncm{\bb}{\bfm{b}}
\ncm{\bg}{\bfm{g}}
\ncm{\be}{\bfm{e}}
\ncm{\bl}{\bfm{l}}
\ncm{\bp}{\bfm{p}}
\ncm{\bq}{\bfm{q}}
\ncm{\br}{\bfm{r}}
\ncm{\bs}{\bfm{s}}
\ncm{\bx}{\bfm{x}}
\ncm{\by}{\bfm{y}}
\ncm{\bz}{\bfm{z}}
\ncm{\balp}{\bfm{\alpha}}
\ncm{\bbe}{\bfm{\beta}}
\ncm{\bxi}{\bfm{\xi}}
\ncm{\bth}{\bfm{\theta}}
\ncm{\bom}{\bfm{\om}}
\ncm{\bmu}{\bfm{\mu}}
\ncm{\bde}{\bfm{\de}}
\ncm{\bva}{\bfm{\va}}
\ncm{\beps}{\bfm{\eps}}
\ncm{\bga}{\bfm{\ga}}
\ncm{\bka}{\bfm{\ka}}
\ncm{\bla}{\bfm{\la}}
\ncm{\bpi}{\bfm{\pi}}
\ncm{\brho}{\bfm{\rho}}
\ncm{\boldeta}{\bfm{\eta}}
\ncm{\bphi}{\bfm{\phi}}
\ncm{\bLa}{\bfm{\Lambda}}
\ncm{\bPi}{\bfm{\Pi}}
\ncm{\bSi}{\bfm{\Si}}
\ncm{\bone}{\bfm{1}}

\ncm{\sbb}{\sbfm{b}}
\ncm{\sbc}{\sbfm{c}}
\ncm{\sbd}{\sbfm{d}}
\ncm{\sbn}{\sbfm{n}}
\ncm{\sbC}{\sbfm{C}}
\ncm{\sbM}{\sbfm{M}}
\ncm{\sbX}{\sbfm{X}}
\ncm{\sbw}{\sbfm{w}}
\ncm{\sbx}{\sbfm{x}}
\ncm{\sby}{\sbfm{y}}
\ncm{\subu}{\scrmath{\ubu}}
\ncm{\subv}{\scrmath{\ubv}}
\ncm{\subw}{\scrmath{\ubw}}
\ncm{\subx}{\scrmath{\ubx}}
\ncm{\suby}{\scrmath{\uby}}
\ncm{\subX}{\scrmath{\ubX}}

\ncm{\hbe}{\hat{\beta}}
\ncm{\heta}{\hat{\eta}}
\ncm{\hth}{\hat{\theta}}
\ncm{\hbth}{\hat{\bth}}
\ncm{\hs}{\hat{s}}
\ncm{\hN}{\hat{N}}
\ncm{\hF}{\hat{F}}
\ncm{\hI}{\hat{I}}
\ncm{\hP}{\hat{P}}
\ncm{\htau}{\hat{\tau}}
\ncm{\hla}{\hat{\lambda}}
\ncm{\hmu}{\hat{\mu}}
\ncm{\hpi}{\hat{\pi}}

\ncm{\mast}{m^\ast}
\ncm{\cast}{c^\ast}
\ncm{\fast}{f^\ast}
\ncm{\siast}{\si^\ast}
\ncm{\psiast}{\psi^\ast}
\ncm{\tsiast}{\tilde{\si}^\ast}
\ncm{\alfast}{\alpha^\ast}
\ncm{\tkaast}{\tilde{\kappa}^\ast}
\ncm{\Xast}{X^\ast}
\ncm{\Yast}{Y^\ast}

\ncm{\ap}{a^\prime}
\ncm{\hp}{h^\prime}
\ncm{\ip}{i^\prime}
\ncm{\jp}{j^\prime}
\ncm{\kp}{k^\prime}
\ncm{\lp}{l^\prime}
\ncm{\np}{n^\prime}
\ncm{\npr}{n^\prime}
\ncm{\qp}{q^\prime}
\ncm{\rp}{r^\prime}
\ncm{\spr}{s^\prime}
\ncm{\up}{u^\prime}
\ncm{\vp}{v^\prime}
\ncm{\wpr}{w^\prime}
\ncm{\xp}{x^\prime}
\ncm{\yp}{y^\prime}
\ncm{\zp}{z^\prime}
\ncm{\Cp}{C^\prime}
\ncm{\Gp}{G^\prime}
\ncm{\Ip}{I^\prime}
\ncm{\Mp}{M^\prime}
\ncm{\Np}{N^\prime}
\ncm{\Npr}{N^\prime}
\ncm{\Tp}{T^\prime}
\ncm{\gap}{\ga^\prime}
\ncm{\phpr}{\phi^\prime}

\ncm{\wbis}{w^{\prime\prime}}

\ncm{\tih}{\tilde{h}}
\ncm{\tZ}{\tilde{Z}}
\ncm{\tA}{\tilde{A}}
\ncm{\tC}{\tilde{C}}
\ncm{\tD}{\tilde{D}}
\ncm{\tF}{\tilde{F}}
\ncm{\tI}{\tilde{I}}
\ncm{\tN}{\tilde{N}}
\ncm{\tY}{\tilde{Y}}
\ncm{\tmu}{\tilde{\mu}}
\ncm{\tOm}{\tilde{\Omega}}
\ncm{\tnu}{\tilde{\nu}}
\ncm{\tsi}{\tilde{\sigma}}
\ncm{\tal}{\tilde{\alpha}}
\ncm{\tbeta}{\tilde{\beta}}
\ncm{\tde}{\tilde{\delta}}
\ncm{\tka}{\tilde{\ka}}
\ncm{\txi}{\tilde{\xi}}
\ncm{\tmathV}{\tilde{\V}}
\ncm{\tV}{\tilde{V}}
\ncm{\tr}{\tilde{r}}
\ncm{\tu}{\tilde{u}}
\ncm{\tw}{\tilde{w}}
\ncm{\twpr}{\tilde{w}^\prime}
\ncm{\tb}{\tilde{b}}
\ncm{\td}{\tilde{d}}
\ncm{\tp}{\tilde{p}}
\ncm{\tf}{\tilde{f}}
\ncm{\tn}{\tilde{n}}
\ncm{\tS}{\tilde{S}}
\ncm{\tL}{\tilde{L}}
\ncm{\tl}{\tilde{l}}
\ncm{\tP}{\tilde{P}}
\ncm{\tSmath}{\tilde{\mathbb{S}}}
\ncm{\tT}{\tilde{T}}
\ncm{\tK}{\tilde{K}}
\ncm{\tva}{\tilde{\va}}
\ncm{\tla}{\tilde{\la}}
\ncm{\tpi}{\tilde{\pi}}
\ncm{\trho}{\tilde{\rho}}
\ncm{\tbom}{\tilde{\bfm{\om}}}
\ncm{\tbxi}{\tilde{\bfm{\xi}}}
\ncm{\tbrho}{\tilde{\bfm{\rho}}}
\ncm{\tbg}{\tilde{\bg}}
\ncm{\tbb}{\tilde{\bb}}
\ncm{\tbr}{\tilde{\br}}
\ncm{\tbf}{\tilde{\bfm{f}}}
\ncm{\tbD}{\tilde{\bfm{D}}}
\ncm{\tbH}{\tilde{\bfm{H}}}
\ncm{\tbone}{\tilde{\bfm{1}}}
\ncm{\tbe}{\tilde{\bfm{e}}}

\ncm{\baW}{\bar{W}}
\ncm{\bacW}{\bar{\cW}}
\ncm{\bacV}{\bar{\cV}}
\ncm{\baf}{\bar{f}}
\ncm{\bah}{\bar{h}}
\ncm{\ban}{\bar{n}}
\ncm{\bap}{\bar{p}}
\ncm{\bav}{\bar{v}}
\ncm{\baw}{\bar{w}}
\ncm{\baZ}{\bar{Z}}
\ncm{\baY}{\bar{Y}}
\ncm{\baS}{\bar{S}}
\ncm{\baH}{\bar{H}}
\ncm{\baA}{\bar{A}}
\ncm{\baD}{\bar{D}}
\ncm{\baC}{\bar{C}}
\ncm{\baN}{\bar{N}}
\ncm{\baQ}{\bar{Q}}
\ncm{\bal}{\bar{l}}
\ncm{\bam}{\bar{m}}
\ncm{\bae}{\bar{e}}
\ncm{\bacR}{\bar{{\cal R}}}
\ncm{\bacP}{\bar{{\cal P}}}
\ncm{\babe}{\bar{\beta}}
\ncm{\baka}{\bar{\kappa}}
\ncm{\bamu}{\bar{\mu}}
\ncm{\banu}{\bar{\nu}}
\ncm{\bade}{\bar{\de}}
\ncm{\bala}{\bar{\la}}
\ncm{\baga}{\bar{\ga}}
\ncm{\barho}{\bar{\rho}}
\ncm{\babf}{\bar{\bfm{f}}}
\ncm{\babD}{\bar{\bfm{D}}}
\ncm{\babA}{\bar{\bfm{A}}}
\ncm{\babQ}{\bar{\bfm{Q}}}
\ncm{\babW}{\bar{\bfm{W}}}
\ncm{\babh}{\bar{\bfm{h}}}
\ncm{\babr}{\bar{\bfm{r}}}
\ncm{\babde}{\bar{\bfm{\de}}}
\ncm{\babrho}{\bar{\bfm{\rho}}}
\ncm{\babone}{\bar{\bfm{1}}}

\ncm{\chnu}{\check{\nu}}

\ncm{\uC}{\underline{C}}
\ncm{\ucH}{\underline{\cH}}
\ncm{\ucX}{\underline{\cX}}
\ncm{\ubx}{\underline{\bx}}
\ncm{\ubu}{\underline{\bu}}
\ncm{\ubv}{\underline{\bv}}
\ncm{\ubw}{\underline{\bw}}
\ncm{\ubX}{\underline{\bX}}
\ncm{\uby}{\underline{\by}}
\ncm{\ubY}{\underline{\bY}}

\ncm{\ocH}{\overline{\cH}}

\ncm{\sca}{\scr{a}}
\ncm{\scn}{\scr{n}}

\ncm{\Lin}{\, \stackrel{\cal L} \in}
\ncm{\Leq}{\, \stackrel{\cal L} =}
\ncm{\Lto}{\, \stackrel{\cal L} \longrightarrow}
\ncm{\pto}{\, \stackrel{p} \longrightarrow}
\ncm{\asto}{\, \stackrel{\rm a.s.} \longrightarrow}
\ncm{\Cov}{\mbox{Cov}}
\ncm{\Var}{\mbox{Var}}
\ncm{\sameord}{\stackrel{\cup}{{\scriptstyle \cap}}}

\ncm{\ith}{i^{\scr{th}}}
\ncm{\jth}{j^{\scr{th}}}
\ncm{\kth}{k^{\scr{th}}}
\ncm{\lth}{l^{\scr{th}}}
\ncm{\Bin}{\mbox{Bin}}
\ncm{\CV}{\mbox{CV}}
\ncm{\Exp}{\mbox{Exp}}
\ncm{\Hyp}{\mbox{Hyp}}
\ncm{\Po}{\mbox{Po}}
\ncm{\mm}{\mbox{mm}}
\ncm{\PD}{\mbox{PD}}
\ncm{\Ctot}{\bar{C}}
\ncm{\Ctottiny}{C_{\mbox{\tiny tot}}}
\ncm{\bzero}{\bfm{0}}
\ncm{\fappr}{\hat{f}}
\ncm{\bappr}{\hat{b}}
\ncm{\laappr}{\hat{\la}}
\ncm{\muappr}{\hat{\mu}}
\ncm{\pappr}{\hat{p}}
\ncm{\piappr}{\hat{\pi}}
\ncm{\kaappr}{\hat{\ka}}
\ncm{\Siappr}{\hat{\Si}}
\ncm{\bSiappr}{\hat{\bSi}}
\ncm{\demax}{\de_{\scr{max}}}
\ncm{\kamax}{\ka_{\scr{max}}}
\ncm{\mumin}{\mu_{\scr{min}}}
\ncm{\hmumin}{\hmu_{\scr{min}}}
\ncm{\Ias}{I_{\scr{as}}}
\ncm{\Ibas}{I_B}
\ncm{\cBall}{\cB_{\scr{all}}}
\ncm{\Inonas}{I_{\scr{nas}}}
\ncm{\Ilong}{I_{\scr{long}}}
\ncm{\Ishort}{I_{\scr{short}}}
\ncm{\cCcoarse}{\cC_{\scr{red}}}
\ncm{\cCmis}{\cC_{\scr{mism}}}
\ncm{\cCmishit}{\cC_{\scr{mismhit}}}
\ncm{\Cmax}{C_{\scr{max}}}

\ncm{\beq}{\begin{equation}}
\ncm{\eeq}{\end{equation}}
\ncm{\beqr}{\begin{eqnarray}}
\ncm{\eeqr}{\end{eqnarray}}
\ncm{\beqrn}{\begin{eqnarray*}}
\ncm{\eeqrn}{\end{eqnarray*}}
\ncm\rthm[1]{\ref{#1}}
\ncm\lb[1]{\label{#1}}
\ncm\re[1]{(\ref{#1})}
\ncm{\slut}{
  {\unskip\nobreak\hfill\penalty100\hskip1em\vadjust{}\nobreak
  \hfill\mbox{$\Box$}\parfillskip=0pt\finalhyphendemerits=0}}

\parindent=0mm
\newcommand*\samethanks[1][\value{footnote}]{\footnotemark[#1]}

\begin{titlepage}

\title{Sharp Lower and Upper Bounds for the Covariance of Bounded Random Variables}

\author{Ola H\"{o}ssjer\footnote{Department of Mathematics, Stockholm University, 106 91 Stockholm, Sweden. Email: ola@math.su.se} \and Arvid Sjölander\footnote{Department of Medical Epidemiology and Biostatistics, Karolinska Institutet, 171 77 Stockholm, Sweden. Email: arvid.sjolander@ki.se}}

\maketitle

\begin{abstract}
\parindent=0pt
In this paper we derive sharp lower and upper bounds for the covariance of two bounded random variables when knowledge about their expected values, variances or both is available. When only the expected values are known, our result can be viewed as an extension of the Bhatia-Davis Inequality for variances. We also provide a number of different ways to standardize covariance. For a binary pair random variables, one of these standardized measures of covariation agrees with a frequently used measure of dependence between genetic variants.  
\par\bigskip
{\bf Keywords:}  Bounded random variables; standardized measure of variation; covariance; lower and upper bounds.  
\end{abstract}

\end{titlepage}

\section{Introduction}

What can be said about the statisical dependency between two random variables $X$ and $Y$, when some information about their marginal distribution is available? The answer to this question depends on the dependency measure being used as well as the type of restrictions that are imposed on the marginal distributions of $X$ and $Y$. The covariance $\Cov(X,Y)=E(XY)-E(X)E(Y)$ is one of the most frequently employed measures of dependence between two random variables $X$ and $Y$, when these are measured on an interval scale. 
The above question can then be phrased as finding lower and upper bounds of $\Cov(X,Y)$ that incorporate any available information about the marginal distributions of $X$ and $Y$. The most well known such covariance bounds 
\beq
-\sqrt{\Var(X)\Var(Y)} \le \Cov(X,Y) \le \sqrt{\Var(X)\Var(Y)}
\lb{CS}
\eeq
follow from the Cauchy-Schwarz Inequality, originally stated by Augustine Louis Cauchy in 1821, and later proved independently by Viktor Bunyakovski and Karl Hermann Schwarz. The lower and upper bounds in \re{CS} only involve the variances $\Var(X)$ and $\Var(Y)$ of $X$ and $Y$, and they are attained when $Y=kX+l$ is a linear function of $X$ with negative and positive slope respectively. A related class of covariance bounds involve not only the marginal distributions of $X$ and $Y$, but more generally the variance of some function $h(X,Y)$ of $X$ and $Y$ (Koop, 1964, Kimeldorf and Sampson, 1973). There is also a large literature on covariance bounds when $X=f(Z)$ and $Y=g(Z)$ are functions of the same random variable $Z$. These results make use of various mathematical tools such as the Hoeffding Inequality (Hoeffding, 1940), Chebyshev's Integral Inequality and Stein operators, see for instance Egozcue (2015), He and Wang (2015), Ernst et al.\ (2019) and references therein.     
\par\smallskip
In this article we consider a pair of bounded random variables $X$ and $Y$, when knowledge about their marginal distributions is given in terms of their expected values $E(X)$ and $E(Y)$, and/or their variances $\Var(X)$ and $\Var(Y)$. Barnett and Dragomir (2004) considered the case when the expected values of $X$ and $Y$ are known, and they derived lower and upper bounds for the covariance of $X$ and $Y$. However, these bounds are not sharp, and may thus include values of the covariance that are logically impossible, given the expected values. We provide sharp lower and upper bounds for $\Cov(X,Y)$ when the expected values of $X$ and $Y$ are known, which extend well known results for binary random variables (Ferguson, 1941, Cureton, 1959, Guilford, 1965, Davenport and El-Sanhurry, 1991). These bounds can also be viewed as a generalization of the Bhatia-Davis Inequality (Bhatia and Davis, 2000), which provides an upper bound on the variance of a bounded random variable, when its expected value is known. We demonstrate that our covariance bounds are attained when the joint distribution of $X$ and $Y$ is discrete, with at most three possible outcomes. We also derive lower and upper bounds of $\Cov(X,Y)$ when the variances of $X$ and $Y$ are known. These bounds are either equal to or truncated versions of the Cauch-Schwarz bounds in \re{CS}, depending on whether the expected values of $X$ and $Y$ are unknown or known. 
\par\smallskip
The covariance bounds that we propose naturally lead to four different standar\-dized measures of covariation between bounded random variables, depending on whether the expected values and variances of these two random variables are known or not. In particular, for binary random variables with known expected values, the corresponding standardized measure of covariation coincides with a measure of dependence used to quantify linkage disequilibrium between two biallelic genetic variants (Lewontin, 1965, Chapter 8 of Thomas, 2004).  
\par\smallskip
Our paper is organized as follows: In Section \ref{Sec:EKnown} we present our new and sharp covariance bounds of $X$ and $Y$ when the expected values but not the variances of these two random variables are known. Then in Section \ref{Sec:VKnown} we derive covariance bounds of $X$ and $Y$ when the variances of these two random variables are known, whereas the expected values are either known or not. The four standardized measures of covariation are introduced in Section \ref{Sec:CovStand}, and finally a discussion in Section \ref{Sec:Disc} concludes.

\section{Covariance bounds when variances are unknown}\lb{Sec:EKnown}

Throughout this article we assume that $a\le X \le b$ and $c\le Y\le d$ are two bounded random variables, restricted by lower and upper bounds $-\infty < a < b < \infty$ and $-\infty< c < d < \infty$ respectively. In this section we will investigate which values are attainable for the covariance $\Cov(X,Y)$ of $X$ and $Y$, when the variances of $X$ and $Y$ are unknown, whereas the expected values $E(X)$ and $E(Y)$ are either known or not. The following theorem treats the case when the expected values are known:   
\par\smallskip
\begin{theorem}\lb{th1} Assume that the expected values $E(X)$ and $E(Y)$ of $a\le X \le b$ and $c\le Y \le d$ are known. Then the covariance of $X$ and $Y$ satisfies 
\beq
\begin{array}{l}
-\min\left[(E(X)-a)(E(Y)-c),(b-E(X))(d-E(Y))\right]\\
\le \Cov(X,Y)\\
\le \min\left[ (E(X)-a)(d-E(Y)),(b-E(X))(E(Y)-c)\right].
\end{array}
\lb{CovBounds}
\eeq
In particular, the lower covariance bound in \re{CovBounds} is attained for a pair $(X,Y)$ of discrete random variables having at most three possible outcomes, with 
\beq
\begin{array}{l}
P(X=x,Y=y) =\\
= \left\{\begin{array}{ll}
\frac{(b-E(X))(d-E(Y))}{(b-a)(d-c)} - \min\left[ \frac{(b-E(X))(d-E(Y))}{(b-a)(d-c)}, \frac{(E(X)-a)(E(Y)-c)}{(b-a)(d-c)} \right], & x=a, y=c,\\
\frac{(E(X)-a)(d-E(Y))}{(b-a)(d-c)} + \min\left[ \frac{(b-E(X))(d-E(Y))}{(b-a)(d-c)}, \frac{(E(X)-a)(E(Y)-c)}{(b-a)(d-c)} \right], & x=b, y=c,\\
\frac{(b-E(X))(E(Y)-c)}{(b-a)(d-c)} + \min\left[ \frac{(b-E(X))(d-E(Y))}{(b-a)(d-c)}, \frac{(E(X)-a)(E(Y)-c)}{(b-a)(d-c)} \right], & x=a, y=d,\\
\frac{(E(X)-a)(E(Y)-c)}{(b-a)(d-c)} - \min\left[ \frac{(b-E(X))(d-E(Y))}{(b-a)(d-c)}, \frac{(E(X)-a)(E(Y)-c)}{(b-a)(d-c)} \right], & x=b, y=d,\\
0, & \mbox{otherwise},
\end{array}\right.
\end{array} \lb{LowBound}
\eeq
whereas the upper covariance bound in \re{CovBounds} is attained for another pair $(X,Y)$ of discrete random variables having at most three possible outcomes, with 
\beq
\begin{array}{l}
P(X=x,Y=y) =\\
= \left\{\begin{array}{ll}
\frac{(b-E(X))(d-E(Y))}{(b-a)(d-c)} + \min\left[ \frac{(E(X)-a)(d-E(Y))}{(b-a)(d-c)}, \frac{(b-E(X))(E(Y)-c)}{(b-a)(d-c)}\right], & x=a, y=c,\\
\frac{(E(X)-a)(d-E(Y))}{(b-a)(d-c)} - \min\left[ \frac{(E(X)-a)(d-E(Y))}{(b-a)(d-c)}, \frac{(b-E(X))(E(Y)-c)}{(b-a)(d-c)}\right], & x=b, y=c,\\
\frac{(b-E(X))(E(Y)-c)}{(b-a)(d-c)} - \min\left[ \frac{(E(X)-a)(d-E(Y))}{(b-a)(d-c)}, \frac{(b-E(X))(E(Y)-c)}{(b-a)(d-c)}\right], & x=a, y=d,\\
\frac{(E(X)-a)(E(Y)-c)}{(b-a)(d-c)} + \min\left[ \frac{(E(X)-a)(d-E(Y))}{(b-a)(d-c)}, \frac{(b-E(X))(E(Y)-c)}{(b-a)(d-c)}\right], & x=b, y=d,\\
0, & \mbox{otherwise}.
\end{array}\right.
\end{array} \lb{UppBound}
\eeq
\end{theorem}
\par\smallskip
It turns out that Theorem \ref{th1} is related to the Bhatia-Davis Inequality for the variance of bounded random variables.  This inequality implies
\beq
\begin{array}{rcl}
\Var(X) &\le & (E(X)-a)(b-E(X)),\\
\Var(Y) &\le & (E(Y)-c)(d-E(Y)).
\lb{BD}
\end{array}
\eeq
Setting $X=Y$, $a=c$, and $b=d$ we find that the upper bound of $\Cov(X,X)=\Var(X)$ in \re{CovBounds} agrees with \re{BD}. 
It is possible to combine the Bhatia-Davis Inequality with the Cauchy-Schwarz Inequality \re{CS}. Indeed, inserting \re{BD} into \re{CS} we deduce 
\beq
\begin{array}{l}
-\sqrt{(E(X)-a)(b-E(X))(E(Y)-c)(d-E(Y))}\\
\le \Cov(X,Y)\\
\le \sqrt{(E(X)-a)(b-E(X))(E(Y)-c)(d-E(Y))}.
\end{array}
\lb{CSBD}
\eeq
It follows from Theorem \ref{th1} that the bounds in \re{CSBD} are at least as wide as those in \re{CovBounds}. We will give precise conditions under which the bounds in \re{CSBD} are strictly wider. To this end, it is helpful to rewrite the expected values of $X$ and $Y$ as
\beq
\begin{array}{rcl}
E(X) &=& a + \alpha (b-a),\\
E(Y) &=& c + \beta (d-c),
\end{array}
\lb{ab}
\eeq
for some constants $0\le \alpha, \beta \le 1$. These numbers quantify the expected values of $X$ and $Y$ on a relative scale, and as the following result shows, they determine when the Cauchy-Schwarz covariance bounds are strictly wider than those of Theorem \ref{th1}:  
\par\smallskip
\begin{corollary}\lb{cor1} Assume $0<\alpha,\beta<1$. The covariance bounds of \re{CSBD} are at least as wide as those of \re{CovBounds}. The lower covariance bound of \re{CSBD} equals the one in \re{CovBounds} if and only if $\alpha+\beta=1$, and then the random vector $(X,Y)$ in \re{LowBound} that attains this lower bound has a two point distribution supported on $(a,d)$ and $(b,c)$. Whenever the lower covariance bounds of \re{CovBounds} and \re{CSBD} differ, the random vector $(X,Y)$ in \re{LowBound} that attains the lower bound of \re{CovBounds} has a three point distribution. The upper covariance bound of \re{CSBD} equals the one in \re{CovBounds} if and only if $\alpha=\beta$, and the random vector $(X,Y)$ in \re{UppBound} that attains this upper bound has a two point distribution supported on $(a,c)$ and $(b,d)$. Whenever the upper bounds of \re{CovBounds} and \re{CSBD} differ,  the random vector $(X,Y)$ in \re{UppBound} that attains the upper bound of \re{CovBounds} has a three point distribution.
\end{corollary}
\par\smallskip
Bartnett and Dragomir (2004) investigated upper and lower bounds of the covariance $\Cov(X,Y)$ of two bounded random variable with known expected values. At the end of Section 8 of their paper, they obtain  
\beq
|\Cov(X,Y) + [b-E(X)][d-E(Y)]| \le (b-a) + (d-c) + (b-a)(d-c).
\lb{BD04}
\eeq
The following result details how \re{BD04} compares to the covariance bounds of Theorem \ref{th1}:
\par\smallskip
\begin{corollary}\lb{cor2} The upper and lower covariance bounds in \re{CovBounds} are strictly sharper than those obtained from \re{BD04}, even when the two terms $b-a$ and $d-c$ are removed from the right-hand side of \re{BD04}.
\end{corollary} 
\par\smallskip
By minimizing (maximizing) the left-hand (right-hand) side of \re{CovBounds} with respect to $E(X)$ and $E(Y)$ , it is possible to derive sharp lower (upper) bounds of $\Cov(X,Y)$ when the expected values are unknown:
\par\smallskip
\begin{corollary}\lb{cor3} Assume that no information is available about the distribution of $X$ and $Y$, except that $a\le X \le b$ and $c\le Y \le d$. Then the covariance of $X$ and $Y$ satisfies 
\beq
-\frac{1}{4}(b-a)(d-c) \le \Cov(X,Y) \le \frac{1}{4}(b-a)(d-c).
\lb{CovBoundsabcd}
\eeq
\end{corollary}

\section{Covariance bounds when variances are known} \lb{Sec:VKnown}

In this section we will assume that the variances $\Var(X)$ and $\Var(Y)$ of $X$ and $Y$ are known. To begin with, we also assume that the expected values $E(X)$ and $E(Y)$ are known. The following result unifies Theorem \ref{th1} with the Cauchy-Schwarz Inequality \re{CS}:  
\par\smallskip
\begin{theorem}\lb{th2} Assume that the expected values $E(X)$ and $E(Y)$, as well as the variances $\Var(X)$ and $\Var(Y)$, of $a\le X \le b$ and $c \le Y \le d$ are known. Then 
\beq
\begin{array}{l}
-\min\left[\sqrt{\Var(X)\Var(Y)},(E(X)-a)(E(Y)-c),(b-E(X))(d-E(Y))\right]\\
\le \Cov(X,Y)\\
\le \min\left[\sqrt{\Var(X)\Var(Y)},(E(X)-a)(d-E(Y)),(b-E(X))(E(Y)-c)\right]
\end{array}
\lb{CovBoundsVar}
\eeq
provide sharp lower and upper bounds for the covariance of $X$ and $Y$.
\end{theorem}
\par\smallskip
\begin{example}[Three point distributions]\lb{ex1} {\rm To gain intuition for the results in Theorems \ref{th1} and \ref{th2}, and their relation to the Cauchy-Schwarz bounds, it is instructive to consider  the special case when both $X$ and $Y$ have three-point distributions, as follows. Let $X\in\{a,E(X),b\}$ and $Y\in\{c,E(Y),d\}$ with $P(X=E(X))=r$, $P(Y=E(Y))=s$, and 
$$
\begin{array}{rcl}
P(X=a|X\in\{a,b\})&=&1-\alpha,\\
P(X=b|X\in\{a,b\})&=&\alpha,\\
P(Y=c|Y\in\{c,d\})&=&1-\beta,\\
P(Y=d|Y\in\{c,d\})&=&\beta.
\end{array}
$$
With these figures we have that the expected values of $X$ and $Y$ are given by \re{ab}, whereas
\beq
\begin{array}{rcl}
\textrm{Var}(X) &=& (1-r)(b-a)^2\alpha(1-\alpha),\\
\textrm{Var}(Y) &=& (1-s)(d-c)^2\beta(1-\beta).
\end{array}
\lb{EVThreePoint}
\eeq 
Define the log odds $\psi_X=\textrm{log}(\alpha/(1-\alpha))$ and $\psi_Y=\textrm{log}(\beta/(1-\beta))$. Computing the ratio between the lower bound for $\textrm{Cov}(X,Y)$ in Theorem \ref{th1} and the Cauchy-Schwarz lower bound, gives, after some algebra
\begin{eqnarray}
\label{eq:corr1}
&&\frac{-\textrm{min}[(E(X)-a)(E(Y )-c), (b-E(X))(d-E(Y ))]}{-\sqrt{\textrm{Var}(X)\textrm{Var}(Y)}}\nonumber\\
&&\phantom{1}=\sqrt{\frac{\textrm{min}\left[\textrm{exp}(\psi_X+\psi_Y),\textrm{exp}(-\psi_X-\psi_Y)\right]}{(1-r)(1-s)}}.
\end{eqnarray}
Similarly, computing the ratio between the upper bound for $\textrm{Cov}(X,Y)$ in Theorem \ref{th1} and the Cauchy-Schwarz upper bound, gives
\begin{eqnarray}
\label{eq:corr2}
&&\frac{\textrm{min}[(E(X)-a)(d-E(Y )), (b-E(X))(E(Y )-c)]}{\sqrt{\textrm{Var}(X)\textrm{Var}(Y)}}\nonumber\\
&&\phantom{1}=\sqrt{\frac{\textrm{min}\left[\textrm{exp}(\psi_X-\psi_Y),\textrm{exp}(-\psi_X+\psi_Y)\right]}{(1-r)(1-s)}}.
\end{eqnarray}
If either $r$ or $s$ approach 1, then the ratios in (\ref{eq:corr1}) and (\ref{eq:corr2}) both approach infinity. Hence, when most of the probability mass is located at the mean for either $X$ or $Y$, the Cauchy-Schwarz bounds tend to be more informative (i.e. narrower) than the bounds in Theorem \ref{th1}, or equivalently, the Cauchy-Schwarz bounds will appear in Theorem \ref{th2}. Conversely, if both $r$ and $s$ approach 0, then the ratios in (\ref{eq:corr1}) and (\ref{eq:corr2}) approach their numerators, respectively, which are both $\leq 1$. Hence, when most of the the probability mass is located at the extreme ends for both $X$ and $Y$, the bounds in Theorem \ref{th1} tend to be more informative than the Cauchy-Schwartz bounds, and therefore the bounds of Theorem \ref{th1} will also appear in Theorem \ref{th2}. An exception from the latter occurs for the lower bound when $\psi_X=-\psi_Y$ (or $\alpha = 1-\beta$), that is, when $X$ and $Y$ have opposite skews (cf.\ Corollary \ref{cor1}). In this case the numerator of (\ref{eq:corr1}) is equal to 1, which implies that the lower bound in Theorem \ref{th1} is never more informative than the Cauchy-Schwartz lower bound. A similar exception occurs for the upper bound when $\psi_X=\psi_Y$ (or $\alpha=\beta$), that is, when $X$ and $Y$ have the same skews (cf.\ Corollary \ref{cor1}). In this case the numerator of (\ref{eq:corr2}) is equal to 1, which implies that the upper bound in Theorem \ref{th1} is never more informative than the Cauchy-Schwartz upper bound.}
\end{example}
\par\smallskip
\begin{example}[Continuous distributions]\lb{ex2} {\rm In order to illustrate the difference between Theorems \ref{th1} and \ref{th2} for continuous random variables, assume that rescaled versions of $X$ and $Y$ have beta distributions. Given numbers $0<\alpha,\beta,r,s < 1$, we postulate  
\beq
\begin{array}{rcl}
(X-a)/(b-a) &\sim & \mbox{Beta}(\alpha r/(1-r),(1-\alpha)r/(1-r)),\\
(Y-c)/(d-c) &\sim & \mbox{Beta}(\beta s/(1-s),(1-\beta)s/(1-s)),
\end{array}
\lb{XYBeta}
\eeq
where the limits $r\to 0, s\to 0$ ($r\to 1, s\to 1$) correspond to the same two point (one point) distributons of $X$ and $Y$ as in Example \ref{ex1}. Using formulas for the expected value and variance of a beta distribution, it follows that the expected values and variances of $X$ and $Y$ are the same as in Example \ref{ex1} (cf.\ \re{ab} and \re{EVThreePoint}), for any values of $\alpha$, $\beta$, $r$, and $s$. Therefore, the ratios between the bounds of Theorem \ref{th1}, and the corresponding Cauchy-Schwarz bounds, are the same as in \re{eq:corr1}-\re{eq:corr2}. }  
\end{example}
\par\smallskip
By minimizing (maximizing) the left-hand (right-hand) side of \re{CovBounds} it is possible to derive lower (upper) bounds of the covariance of $X$ and $Y$ when the variances but not the expected values of these two random variables are known: 
\par\smallskip 
\begin{corollary}\lb{cor4} Assume that the variances $\Var(X)$ and $\Var(Y)$ of $a\le X\le b$ and $c\le Y \le d$ are known. Then the Cauchy-Schwarz Inequality \re{CS} provide sharp lower and upper bounds for the covariance of $X$ and $Y$. 
\end{corollary}

\section{Standardized Measures of Covariation}\lb{Sec:CovStand}

In this section we will present four different ways of standardizing the covariance of $a\le X \le b$ and $c\le Y\le d$, so that all values in $[-1,1]$ are possible for the standardized measure. The form of these standardized covariances will depend on whether the expected values and variances of $X$ and $Y$ are known or not.  

\subsection{No moments known}

When neither the expected values nor the variances of $X$ and $Y$ are known we use Corollary \ref{cor3} and introduce
$$
D(X,Y) = \frac{C(X,Y)}{(b-a)(d-c)/4}.
$$

\subsection{Variances known}

When the variances but not the expected values of $X$ and $Y$ are known, we employ Corollary \ref{cor4} and use the ordinary correlation coefficient 
$$
r(X,Y) = \frac{\Cov(X,Y)}{\sqrt{\Var(X)\Var(Y)}}
$$
as a standardized version of the covariance. 

\subsection {Exected values known}

Assume that the expected values but not variances of $X$ and $Y$ are known. Then we use Theorem \ref{th1} and define
\beq
D^\prime(X,Y) = \left\{\begin{array}{ll}
\frac{\scr{Cov}(X,Y)}{\min\left[(E(X)-a)(E(Y)-c),(b-E(X))(d-E(Y))\right]}, & \Cov(X,Y) < 0, \\
\frac{\scr{Cov}(X,Y)}{\min\left[ (E(X)-a)(d-E(Y)),(b-E(X))(E(Y)-c)\right]}, & \Cov(X,Y) \ge 0
\end{array}\right.
\lb{Dprime}
\eeq
as a standardized covariance. In particular, when $X$ and $Y$ have two point distributions on $\{a,b\}$ and $\{c,d\}$, $D^\prime(X,Y)$ is a well known measure of dependence (Ferguson, 1941). In genetic epidemiology it is a freqeuently used measure of linkage disequilibrium betweeen two biallelic genetic variants (Lewontin, 1965). 

\subsection{Expected values and variances known}

If the expected values and variances of $X$ of $Y$ are known, it is natural to use Theorem \ref{th2} for standardizing the covariance of $X$ and $Y$. This amounts to a definition  
\beq
D^{\prime\prime}(X,Y) = \left\{\begin{array}{ll}
\frac{\scr{Cov}(X,Y)}{\min\left[\sqrt{\scr{Var}(X)\scr{Var}(Y)},(E(X)-a)(E(Y)-c),(b-E(X))(d-E(Y))\right]}, & \Cov(X,Y) < 0, \\
\frac{\scr{Cov}(X,Y)}{\min\left[\sqrt{\scr{Var}(X)\scr{Var}(Y)},(E(X)-a)(d-E(Y)),(b-E(X))(E(Y)-c)\right]}, & \Cov(X,Y) \ge 0.
\end{array}\right.
\lb{D}
\eeq

\subsection{Relations between the standardized measures of covariation}

Our four measures of standardized covariation have a partial ordering
\beq
\begin{array}{rcl}
|D(X,Y)| &\le & \min(|r(X,Y)|,|D^\prime(X,Y)|),\\
|D^{\prime\prime}(X,Y)| &\ge & \max(|r(X,Y)|,|D^\prime(X,Y)|).
\end{array}
\lb{Ordering}
\eeq
There is however no general ordering between $|r(X,Y)|$ and $|D^\prime(X,Y)|$. Although $|r(X,Y)|\le |D^\prime(X,Y)|$ holds for binary random variables, we recall from Examples \ref{ex1}-\ref{ex2} that this inequality sometimes goes in the other direction when $X$ and $Y$ have three point distributions or beta distributions.  

\section{Discussion}\lb{Sec:Disc}

In this paper we derived sharp lower and upper bounds for the covariance of two bounded random variables $X$ and $Y$ when their expected values and/or their variances, are known. This resulted in various ways of standardizing covariances, some of which are well known, whereas others are new. 
A number of extensions are of interest. A first extension is to find the minimum and maximum covariance of two bounded random variables under other moment constraints than expected values and variances. More generally, it would be of interest to derive covariance bounds under various types of restrictions on the marginal distributions of $X$ and $Y$.   
A second extension is to obtain bounds for other types of dependency measures between $X$ and $Y$, under various restrictions on the marginal distributions of these two random variables. Examples of alternative dependency measures include the kappa statistic (Cohen, 1960) and proportional reduction in entropy (Theil, 1970) for nominal random variables and the gamma statistic (Goodman and Kruskal, 1954) for ordinal random variables. 

\setcounter{equation}{0}
\renewcommand{\thesection}{A}
\renewcommand{\theequation}{A.\arabic{equation}}
\section{Appendix}

\subsection{Proofs from Section \ref{Sec:EKnown}.}

\begin{proof}[Proof of Theorem \ref{th1}]. Since the covariance operator as well as the lower and upper bounds of \re{CovBounds} are bilinear, equation \re{CovBounds} is invariant with respect to linear transformations of $X$ and $Y$. We may therefore without loss of generality assume $a=c=0$ and $b=d=1$. Thus our objective is to prove 
\beq
\begin{array}{l}
-\min\left[E(X)E(Y),(1-E(X))(1-E(Y))\right]\\
\le \Cov(X,Y)\\
\le \min\left[ E(X)(1-E(Y)),(1-E(X))E(Y)\right],
\end{array}
\lb{CovBounds2}
\eeq
for pairs $(X,Y)$ of random variables satsifying $0\le X,Y \le 1$, or equivalently
\beq
\begin{array}{l}
E(X)E(Y)-\min\left[E(X)E(Y),(1-E(X))(1-E(Y))\right]\\
\le E(XY)\\
\le E(X)E(Y)+\min\left[ E(X)(1-E(Y)),(1-E(X))E(Y)\right].
\end{array}
\lb{EXYBounds}
\eeq
Moreover, we also need to show that the lower (upper) bounds of \re{EXYBounds} are attained by a binary pair of random variables satisfying \re{LowBound} and \re{UppBound} respectively, with $a=c=0$ and $b=d=1$. 
\par\smallskip
Given any pair $0\le X,Y \le 1$ of random variables, there is a corresponding binary pair of random variables $\Xast,\Yast\in \{0,1\}$ satisfying
$$
\begin{array}{rcl} 
P(\Xast=0,\Yast=0) &=& 1 - E(X) - E(Y) + E(XY),\\
P(\Xast=1,\Yast=0) &=& E(X) - E(XY),\\
P(\Xast=0,\Yast=1) &=& E(Y) - E(XY),\\
P(\Xast=1,\Yast=1) &=& E(XY),
\end{array}
$$
and consequently
\beq
\begin{array}{rcl}
E(\Xast) &=& E(X),\\
E(\Yast) &=& E(Y),\\
E(\Xast\Yast) &=& E(XY).
\end{array}
\lb{1}
\eeq
In view of \re{1}, maximizing (minimizing) $E(XY)$ over all pairs $(X,Y)$ of random variables with fixed expected values is equivalent to maximizing (minimizing) $E(\Xast \Yast)$ over all binary pairs $(\Xast,\Yast)$ of random variables with $E(\Xast)=E(X)$ and $E(\Yast)=E(Y)$. It therefore suffices to establish \re{EXYBounds}, with $E(XY)$ replaced by $E(\Xast \Yast)$, for all binary pairs $\Xast,\Yast\in\{0,1\}$ of random variables with $E(\Xast)=E(X)$ and $E(\Yast)=E(Y)$. In order to simplify notation we introduce 
$$
\begin{array}{rcl}
p_{00} &=& P(\Xast=0,\Yast=0),\\
p_{10} &=& P(\Xast=1,\Yast=0),\\
p_{01} &=& P(\Xast=0,\Yast=1),\\
p_{11} &=& P(\Xast=1,\Yast=1).
\end{array}
$$
We want to find those column vectors $\bp=(p_{00},p_{10},p_{01},p_{11})^T$ that maximize (minimize) $E(\Xast \Yast) = p_{11}$ subject to the constraints 
\beq
\begin{array}{rcl}
p_{10}+p_{11} &=& E(X),\\
p_{01}+p_{11} &=& E(Y),\\
p_{00}+p_{10}+p_{01}+p_{11} &=& 1,\\
p_i &\ge& 0, \,\,\,\,\, i\in \{00,10,01,11\}.
\end{array}
\lb{Constr}
\eeq
The first three linear constraints of \re{Constr} can be written in marix form as
\beq
\bA\bp = \left(\begin{array}{c} E(X) \\ E(Y) \\ 1 \end{array}\right), 
\lb{LinConstr}
\eeq
where 
$$
\bA = \left(\begin{array}{cccc}
0 & 1 & 0 & 1 \\
0 & 0 & 1 & 1 \\
1 & 1 & 1 & 1 
\end{array}\right).
$$
Define a column vector $\bq = (q_{00},q_{10},q_{01},q_{11})^T$ that corresponds to a pair of independent random variables $\Xast$ and $\Yast$, so that
$$
\begin{array}{rcl}
q_{00} &=& P(\Xast=0,\Yast=0) = (1-E(X))(1-E(Y)),\\
q_{10} &=& P(\Xast=1,\Yast=0) = E(X)(1-E(Y)),\\
q_{01} &=& P(\Xast=0,\Yast=1) = (1-E(X))E(Y),\\
q_{11} &=& P(\Xast=1,\Yast=1) = E(X)E(Y).
\end{array}
$$
It is clear that $\bq$ satisfies the constraints imposed in \re{LinConstr}. Therefore,  any solution of the linear inhomogeneous equation \re{LinConstr} is of the form $\bp=\bq+\bv$, where $\bv=(v_{00},v_{10},v_{01},v_{11})^T$ is any solution of the corresponding homogeneous linear equation $\bA\bv=\bzero$, i.e.
$$
\begin{array}{rcl}
v_{10} + v_{11} &=& 0,\\
v_{01} + v_{11} &=& 0,\\
v_{00} + v_{10} + v_{01} + v_{11} &=& 1. 
\end{array}
$$
This is equivalent to $\bv = c(1,-1,-1,1)$ for some constant $c\in\R$. Thus we want to maximize (minimize) $p_{11}=q_{11} + c$ subject to the inequality constraints of \re{Constr}, i.e.
\beq
\begin{array}{rcl}
q_{00} + c &\ge & 0,\\
q_{10} - c &\ge & 0,\\
q_{01} - c &\ge & 0,\\
q_{11} + c &\ge & 0.
\end{array}
\lb{cConstr}
\eeq
It is clear that $q_{11}+c$ is minimized when $c$ is chosen as small as possible, and yet satisfies \re{cConstr}. This corresponds to $c=-\min(q_{00},q_{11})$ and 
$$
\begin{array}{rcl}
E(\Xast \Yast)_{\scr{min}} &=& q_{11} - \min(q_{00},q_{11})\\
&=& E(XY) - \min\left[ (1-E(X))(1-E(Y)),E(X)E(Y)\right],
\end{array}
$$
in agreement with the lower bound of \re{EXYBounds}. Analogously, $q_{11}+c$ is maximized for $c=\min(q_{10},q_{01})$, corresponding to  
$$
\begin{array}{rcl}
E(\Xast \Yast)_{\scr{max}} &=& q_{11} + \min(q_{10},q_{01})\\
&=& E(XY) + \min\left[ E(X)(1-E(Y)),(1-E(X))E(Y)\right],
\end{array}
$$
in agreement with the upper bound of \re{EXYBounds}. The proof is finalized by noticing that the two vectors $\bp=\bq-\min(q_{00},q_{11})\bv$ and $\bp=\bq+\min(q_{10},q_{01})\bv$ correspond to the bivariate distributions of $(X,Y)$ in \re{LowBound} and \re{UppBound} respectively, when $a=c=0$ and $b=d=1$. 
\end{proof}

\begin{proof}[Proof of Corollary \ref{cor1}.] It is helpful to rewrite \re{CovBounds} and \re{CSBD} as 
\beq
\begin{array}{l}
-(b-a)(d-c)\min\left[\alpha\beta,(1-\alpha)(1-\beta)\right]\\
\le \Cov(X,Y)\\
\le (b-a)(d-c)\min\left[ \alpha (1-\beta),(1-\alpha)\beta\right]
\end{array}
\lb{CovBoundsStand}
\eeq
and
\beq
\begin{array}{l}
-(b-a)(d-c)\sqrt{\alpha (1-\alpha)\beta (1-\beta)}\\
\le \Cov(X,Y)\\
\le (b-a)(d-c)\sqrt{\alpha (1-\alpha)\beta (1-\beta)}
\end{array}
\lb{CSBDStand}
\eeq
respectively. The lower bound of \re{CSBDStand} is at least as small as that in \re{CovBoundsStand}, since
$$
-\sqrt{\alpha (1-\alpha)\beta (1-\beta)} = - \sqrt{\alpha\beta \cdot (1-\alpha)(1-\beta)} \le - \min\left[\alpha\beta,(1-\alpha)(1-\beta)\right].
$$
Moreover, it is clear that the lower bounds of \re{CovBoundsStand} and \re{CSBDStand} agree if and only if 
\beq
\alpha\beta=(1-\alpha)(1-\beta) \Longleftrightarrow \frac{\alpha}{1-\alpha} = \frac{(1-\beta)}{\beta}. 
\lb{alphabeta}
\eeq
Since $x\to f(x)=x/(1-x)$ is strictly increasing on $(0,1)$ it follows that \re{alphabeta} is equivalent to $\alpha=1-\beta$. Moreover, when $\alpha=1-\beta$, the random vector $(X,Y)$ of \re{LowBound} has a two point distribution, since $P(X=x,Y=y)=0$ when $(x,y)$ equals $(a,c)$ and $(b,d)$. This concludes the proof for the lower covariance bounds \re{CovBounds} and \re{CSBD}. The proof for the upper covariance bounds is analogous.  
\end{proof}

\begin{proof}[Proof of Corollary \ref{cor2}.] Denote the upper and lower covariance bounds of \re{CovBounds} by $U$ and $L$, whereas those in \re{BD04} are denoted $U_{\scr{BD}}$ and $L_{\scr{BD}}$ respectively.  We will start comparing the two upper covariance bounds. Recall from \re{CovBoundsStand} that 
$$
U = \min\left[ \alpha (1-\beta),(1-\alpha)\beta\right](b-a)(d-c),
$$
whereas the upper covariance bound of \re{BD04} takes the form
$$
\begin{array}{rcl}
U_{\scr{BD}} &=& - [b-E(X)][d-E(Y)]| + (b-a) + (d-c) + (b-a)(d-c)\\
&=& (b-a) + (d-c) + (\alpha + \beta - \alpha\beta) (b-a)(d-c),
\end{array}
$$
where in the last step we made use of \re{ab}. Hence
$$
\begin{array}{rcl}
U_{\scr{BD}} - U &=&  (b-a) + (d-c) +  \left[\alpha + \beta - \min(\alpha,\beta)\right] (b-a)(d-c)\\
&>& (b-a) + (d-c)\\
& >&  0.
\end{array}
$$
For the lower covariance bounds we similarly derive 
$$
\begin{array}{rcl}
L &=& - \min\left[\alpha\beta,(1-\alpha)(1-\beta)\right](b-a)(d-c),\\
L_{\scr{BD}} &=& - (b-a) - (d-c) - \left[1+(1-\alpha)(1-\beta)\right](b-a)(c-d).
\end{array}
$$
Consequently 
$$
\begin{array}{rcl}
L_{\scr{BD}} - L &=& - (b-a) - (d-c)  \\
&-& \left\{1 +(1-\alpha)(1-\beta) - \min\left[\alpha\beta,(1-\alpha)(1-\beta)\right]\right\} (b-a)(d-c)\\
&<& - (b-a) - (d-c) - (b-a)(d-c)\\
&<& - (b-a) - (d-c)\\
&<& 0.
\end{array} 
$$
\end{proof}

\begin{proof}[Proof of Corollary \ref{cor3}.] We will only verify the upper bound of \re{CovBoundsabcd}, since the proof of the lower bound is analogous. By maximizing the upper bound in \re{CovBoundsabcd} with respect to $E(X)$ and $E(Y)$, and making use of the parametrization \re{ab}, it follows that 
\beq
\begin{array}{rcl}
\Cov(X,Y) &\le & (b-a)(d-c) \sup_{\alpha,\beta} \min\left[ \alpha (1-\beta), (1-\alpha)\beta\right]\\
&\le & (b-a)(d-c)  \sup_{\alpha} \sqrt{\alpha (1-\alpha)} \cdot \sup_\beta \sqrt{\beta(1-\beta)}\\
&=& (b-a)(d-c) \sqrt{ 0.5(1-0.5)} \cdot \sqrt{0.5(1-0.5)}\\
&=& \frac{1}{4} (b-a)(d-c),
\end{array}
\lb{CovXY}
\eeq
where in the second step we invoked Corollary \ref{cor1}. Thus the upper bound of $\Cov(X,Y)$ is at most equal to $(b-a)(d-c)/4$. The fact that the upper bound of $\Cov(X,Y)$ indeed has this value follows from the fact that there is equality in the second step of \re{CovXY} when $\alpha=\beta=0.5$, the values of $\alpha$ and $\beta$ for which the maximum in the third step of \re{CovXY} was attained.  
\end{proof}

\subsection{Proofs from Section \ref{Sec:VKnown}}

\begin{proof}[Proof of Theorem \ref{th2}.] As in the proof of Theorem \ref{th1} we assume without loss of generality that $a=c=0$ and $b=d=1$. Hence we need to prove that 
\beq
\begin{array}{l}
-\min\left[\sqrt{\Var(X)\Var(Y)},E(X)E(Y),(1-E(X))(1-E(Y))\right]\\
\le \Cov(X,Y)\\
\le \min\left[ \sqrt{\Var(X)\Var(Y)},E(X)(1-E(Y)),(1-E(X))E(Y)\right].
\end{array}
\lb{CovBounds3}
\eeq
Once the upper bound of \re{CovBounds3} is established, the lower bound follows from the substitution $X \leftarrow  1-X$, since $\Cov(Y,1-X)=-\Cov(X,Y)$, $E(1-X)=1-E(X)$, and $\Var(1-X)=\Var(X)$. 
\par\smallskip
Hence it sufficies to establish the upper covariance bound 
\beq
\kamax = \min\left[ \sqrt{\Var(X)\Var(Y)},E(X)(1-E(Y)),(1-E(X))E(Y)\right]
\lb{Cmax0}
\eeq
of \re{CovBounds3} and prove that it is sharp. Let
$$
l(x;\ka) = E(Y) + \ka(x-E(X))/\Var(X)
$$
be a line whose interecept and slope involve $E(Y)$, $E(X)$ and $\Var(X)$. Notice in particular that the function $l(x;\ka)$ is constructed in such a way that the random variable $l(X;\ka)$ has expected value $E(l(X;\ka))=E(Y)$ and covariance $\ka=\Cov(X,l(X;\ka))=\{\Var(X)\Var(l(X;\ka))\}^{1/2}$ with $X$. We will first find the number $\tka$ such that $\tY=l(X;\tka)$ with $0\le \tY \le 1$ makes $\Cov(X,\tY)=\tka=\{\Var(X)\Var(\tY)\}^{1/2}$ as large as possible. Note that $E(\tY)=E(Y)$, whereas typically $\Var(\tY)=\Var(Y)$ does not hold. Since $0\le \tY \le 1$ it follows that $l(x;\tka) \in [0,1]$ for all $x$ that belong to the support of the distribution of $X$. Therefore, it is clear that in order to find $\tka$ the support of $X$ should be chosen as small as possible given the pre-specified values of $E(X)$ and $\Var(X)$. A minimal support for the distribution of $X$ is obtained for a two point distribution supported at $0\le x_1 < E(X) < x_2 \le 1$, with $P(X=x_1)=1-p$ and $P(X=x_2)=p$ for some $0\le p \le 1$. In view of the restrictions on the first two moments of $X$ we have that 
\beq
\begin{array}{rcl}
E(X) &=& (1-p)x_1 + p x_2,\\
\Var(X) &=& (1-p)(x_1-E(X))^2 + p(x_2-E(X))^2.
\end{array}
\lb{SE}
\eeq
Moreover, since the slope $\tka>0$ of $l(x;\tka)$ is maximal among all linear functions that map $\{x_1,x_2\}$ to a two point subset of $[0,1]$, it is clear that 
$$
\begin{array}{rcl}
l(x_1;\tka) &=& E(Y) + \tka (x_1-E(X))/\Var(X) \ge 0,\\
l(x_2;\tka) &=& E(Y) + \tka (x_2-E(X))/\Var(X) \le 1,
\end{array}
$$
where at least one of these two inequalities can be replaced by an equality. Consequently, by varying the two point distribution of $X$ we find that 
\beq
\tka = \Var(X) \max_p \min\left(\frac{E(Y)}{E(X)-x_1},\frac{1-E(Y)}{x_2-E(X)}\right),
\lb{tC}
\eeq
since $x_1=x_1(p)$ and $x_2=x_2(p)$ are uniquely determined by $p$ through the system of equations \re{SE}, and therefore the two point distribution of $X$ has only one degree of freedom. 
\par\smallskip
Having defined $\tka$, our next objective is to prove that the upper covariance bound of \re{Cmax0} is given by
\beq
\kamax = \min(\sqrt{\Var(X)\Var(Y)},\tka).
\lb{Cmax}
\eeq
Indeed, if $\ka=\{\Var(X)\Var(Y)\}^{1/2}\le \tka$, because of \re{tC} we can find a random variable $Y=l(X;\ka)$ with $0\le Y \le 1$ and the pre-specified values of $E(Y)$ and $\Var(Y)$, such that $\Cov(X,Y)=\ka$ is maximal. Thus we have found a random variable $Y$ attaining the upper covariance bound $\kamax$ in \re{Cmax} when $\ka=\{\Var(X)\Var(Y)\}^{1/2}\le \tka$, proving that \re{Cmax} is a sharp upper bound of the covariance in this case.  
\par\smallskip
In order to verify \re{Cmax} when $\ka=\{\Var(X)\Var(Y)\}^{1/2}>\tka$ we need to show that $\kamax=\tka$. This follows from the fact that 
\beq
\begin{array}{rcl}
\Cov(X,Y) &=& E(XY) - E(X)E(Y)\\
&=& E(Xm(X)) - E(X)E(Y)\\
& = & E(Xm^\ast(X^\ast)) - E(X)E(Y)\\
& = & E(X^\ast m^\ast(X^\ast)) - E(X)E(Y)\\
& = & \Cov(X^\ast,m^\ast(X^\ast))\\
& \le & \tka,
\end{array}
\lb{Ineq}
\eeq      
where in the second step we introduced $m(x) = E(Y|X=x)$, and in the third step we defined another function $m^\ast(x)$ that attains the two values 
$$
m^\ast(x) = \left\{\begin{array}{ll}
m_1 = \frac{E(Xm(X) | X < E(X))}{E(E(X)-X|X<E(X))}, & x< E(X),\\
m_2 = \frac{E(Xm(X) | X \ge E(X))}{E(X-E(X)|X \ge E(X))}, & x\ge E(X).
\end{array}\right. 
$$
In the fourth step of \re{Ineq} we introduced the binary random variable $X^\ast$, with $P(X^\ast=x_1)=1-p$ and $P(X^\ast=x_2)=p$ for some $0\le x_1 < E(X) < x_2\le 1$ and $p$ that satisfy \re{SE}, so that $E(X^\ast)=E(X)$ and $\Var(\Xast)=\Var(X)$. The inequality in the last step of \re{Ineq} follows from the definition of $\tka$, since $Y^\ast=m^\ast(X^\ast) = l(X^\ast;\ka^\ast)$ satisfies $\ka^\ast = \Cov(X^\ast,Y^\ast) \le \tka$. This concludes the proof of \re{Cmax} when  $\ka=\{\Var(X)\Var(Y)\}^{1/2}>\tka$. 
\par\smallskip
In order to verify that \re{Cmax} is sharp when $\ka=\{\Var(X)\Var(Y)\}^{1/2}>\tka$ we need to find a pair of random variables $X$ and $Y$ with the prescribed expected values and variances, that satisfy $\Cov(X,Y)=\kamax=\tka$. It is possible to choose $X$ as a two-point distribution supported on $x_1$ and $x_2$, with values of $p$, $x_1$, and $x_2$ determined by \re{tC}, and $Y = l(X;\tka) + \va$. The term $\va$, which quantifies a departure from a linear relation between $X$ and $Y$, satisfies $E(\va|X=x_1)=E(\va|X=x_2)=0$ and $(1-p)\Var(\va|X=x_1)+p\Var(\va|X=x_2) = \ka-\tka$.  
\par\smallskip
It remains to verify that \re{Cmax} equals \re{Cmax0}, and this requires an explicit formula for $\tka$. To this end we first note that the upper equation of \re{SE} implies  
$$
p = \frac{E(X)-x_1}{x_2-x_1}.
$$
Insertion of this expression for $p$ into the lower equation of \re{SE} yields
$$
\begin{array}{rcl}
\Var(X) &=& (E(X)-x_1)(x_2-E(X))\\
&=& p(1-p)(x_2-x_1)^2.
\end{array}
$$
As substitution of the last two displayed equations into \re{tC} gives
\beq
\begin{array}{rcl}
\tka &=& \Var(X)\max_p \min\left(\frac{E(Y)}{p(x_2-x_1)},\frac{1-E(Y)}{(1-p)(x_2-x_1)}\right)\\
&=& \sqrt{\Var(X)}\max_p \min\left(\sqrt{\frac{1-p}{p}} E(Y),\sqrt{\frac{p}{1-p}} (1-E(Y)) \right)\\
&=& \sqrt{\Var(X)}\max_\ga \min\left(\ga E(Y), \frac{1-E(Y)}{\ga}\right),
\end{array}
\lb{tCb}
\eeq
where in the last step we introduced $\ga =\{p/(1-p)\}^{1/2}$. The maximization in \re{tCb} ranges over all $p$ (or $\ga$) such that
$$
\begin{array}{rcl}
x_1 &=& E(X) - p(x_2-x_1) = E(X) - \ga\Var(X) \ge 0,\\
x_2 &=& E(X) + (1-p)(x_2-x_1) = E(X) + \ga^{-1}\Var(X) \le 1,
\end{array}
$$ 
which is equivalent to
\beq
\frac{\sqrt{\Var(X)}}{E(X)} \le \ga \le \frac{1-E(X)}{\sqrt{\Var(X)}}.
\lb{ga}
\eeq
In order to further simplifiy \re{tCb}, note that $\ga E(Y) = (1-E(Y))/\ga$ when $\ga = \ga_0=\{(1-E(Y))/E(Y)\}^{1/2}$. We therefore distinguish between three cases, depending on whether $\ga_0$ is located to the left of, to the right of or within the interval \re{ga}. For Case 1 ($\ga_0 \le \Var(X)^{1/2}/E(X)$) we have that
\beq
\begin{array}{rcl}
\tka &=& \sqrt{\Var(X)} \cdot (1-E(Y))/(\sqrt{\Var(X)}/E(X))\\
&=& E(X)(1-E(Y))\\
&=& \min\left(E(X)(1-E(Y)),(1-E(X))E(Y)\right),
\end{array}
\lb{tC1}
\eeq
where in the last step we used
$$
(1-E(X))E(Y) \ge \Var(X)E(Y)/E(X) \ge E(X)(1-E(Y)).
$$
Case 2 ($\ga_0 \ge (1-E(X))/\Var(X)^{1/2}$) is analogous, with  
\beq
\begin{array}{rcl}
\tka &=& \sqrt{\Var(X)} \cdot E(Y) (1-E(X))/\sqrt{\Var(X)}\\
&=& E(Y)(1-E(X))\\
&=& \min\left(E(Y)(1-E(X)),E(X)(1-E(Y))\right),
\end{array}
\lb{tC2}
\eeq 
using
$$
E(X)(1-E(Y)) \ge \Var(X)(1-E(Y))/(1-E(X)) \ge E(Y)(1-E(X))
$$
in the last step. By a similar calculation for Case 3 ($\Var(X)^{1/2}/E(X) \le \ga_0 \le (1-E(X))/\Var(X)^{1/2}$) we find that
\beq
\begin{array}{rcl}
\tka &=& \sqrt{\Var(X)} \cdot \ga_0 E(Y)\\
&=& \sqrt{\Var(X) E(Y) (1-E(Y))}\\
&=& \min\left(\sqrt{\Var(X) E(Y) (1-E(Y))},E(X)(1-E(Y)),(1-E(X))E(Y)\right).
\end{array}
\lb{tC3}
\eeq
Combining \re{tC1}-\re{tC3} with \re{Cmax} we end up with \re{Cmax0} for either of Case 1, Case 2, and Case 3. This completes the proof of the theorem. 
\end{proof}

\begin{proof}[Proof of Corollary \ref{cor4}.] Only the upper bound  $[\Var(X)\Var(Y)]^{1/2}$ of $\Cov(X,Y)$ in Corollary \ref{cor4} will be verified (the lower bound $-[\Var(X)\Var(Y)]^{1/2}$ is derived analogously). In order to verify the upper covariance bound we maximize the upper bound of \re{CovBounds} in Theorem \ref{th2} with respect to $E(X)$ and $E(Y)$ and follow the same line of reasoning as in \re{CovXY}. This gives 
\beq
\begin{array}{rcl}
\Cov(X,Y) & \le & \max_{\alpha,\beta} \min\left[ \sqrt{\Var(X)\Var(Y)}, (b-a)(d-c)\alpha (1-\beta), \right.\\
&& \,\,\,\,\,\,\,\,\,\,\,\,\,\,\, \,\,\,\,\,\,\,\,\,\,\,\,\,\,\,\,\,\,\,\, \left.(b-a)(d-c))(1-\alpha)\beta\right] \\
&=& \min\left[\sqrt{\Var(X)\Var(Y)},(b-a)(d-c)/4\right]\\
&=& \sqrt{\Var(X)\Var(Y)}.
\end{array}
\lb{CovXYCS}
\eeq
The last step of \re{CovXYCS} follows by maximizing the right-hand sides of the Bhatia-Davies Inequalities \re{BD}, with respect to $E(X)$ and $E(Y)$, which gives $\Var(X)\le (b-a)^2/4$ and $\Var(Y)\le (d-c)^2/4$. 
\end{proof}


\begin{thebibliography}{99}

\bibitem{BaDr04} Barnett, N.S.\ and Dragomir, S.S. (2004). Some further inequalities for univariate moments and some new ones for the covariance. \textit{Computers and Mathematics with Applications} \textbf{47}, 23-36. 

\bibitem{BhDa00} Bhatia, R.\ and Davis, C. (2000). A better bound on the variance. \textit{The American Mathematical Monthly} \textbf{107}(4), 353-357.  

\bibitem{Co60} Cohen, J. (1960). A coefficient of agreement for nominal scales. \textit{Educational and Psychological Measurement} \textbf{20}, 37-46.

\bibitem{Cu59} Cureton, E.E. (1959). Note on $\phi/\phi_{\scr{max}}$. \textit{Psychometrika} \textbf{24}, 89-91. 

\bibitem{DaEl91} Davenport, E.C.\ and El-Sanhurry, N.A. (1991). Phi/Phimax: Review and Synthesis. \textit{Educational and Psychological Measurement} \textbf{51}, 821-828.  

\bibitem{Eg15} Egozcue, M. (2015). Some covariance inequalities for non-monotonic functions with applications to mean-variance indifference curves and bank hedging. \textit{Cogent Mathematics} \textbf{2}:991082. 

\bibitem{Er19} Ernst, M., Reinert, G.\ and Swan, Y. (2019). First order covariance inequalitites via Stein's method. arXiv 1906.08372v1. 

\bibitem{Fe40} Ferguson, G.A. (1941). The factorial interpretation of test difficulty. \textit{Psychometrika} \textbf{6}, 323-333.  

\bibitem{GoKr54} Goodman, L.A.\ and Kruskal, W.H. (1954). Measures of association for cross classifications. \textit{Journal of the American Statistical Association} \textbf{49}, 732-764. 

\bibitem{Gu65} Guilford, J.P. (1965). The minimal phi coefficient and the maximal phi. \textit{Educational and Psychological Measurement} \textbf{25},3-8.   

\bibitem{HeWa05} He, Z.\ and Wang, M. (2015). An inequality for covariance with applications. \textit{Journal of Inequalities and Applications} \textbf{2015}:413.   

\bibitem{Ho41} Höffding, W. (1940). Masstabinvariante Korrelationstheorie. \textit{Schriften des Mathematischen Instituts und Instituts for Angewandte Mathematik der Universität Berlin}, \textbf{5}: 181-233.  

\bibitem{KiSa73} Kimeldorf, G.\ and Sampson, A. (1973). A class of covariance inequalities. \textit{Journal of the American Statistical Association} \textbf{68}(341), 228-230. 

\bibitem{Ko64} Koop, J.C. (1964). Some properties of random variables. \textit{Nature} \textbf{203}, 1097-1098. 

\bibitem{Le65} Lewontin, R.C. (1964). The interaction of selection and linkage. I. General considerations; heterotic models. \textit{Genetics} \textbf{49}(1), 49-67. 


\bibitem{Th04} Theil, H. (1970). On the estimation of relationships involving qualitative variables. \textit{American Journal of Sociology} \textbf{76}, 103-154. 

\bibitem{Th04} Thomas, D. (2004). \textit{Statistical Methods in Genetic Epidemiology}. Oxford University Press, Oxford.  

\end{thebibliography}
\end{document}